\documentclass[a4paper,11pt]{amsart}
\usepackage{graphicx}
\usepackage{mathptmx}
\usepackage{amsmath}
\usepackage{amssymb}
\usepackage{enumitem}
\usepackage{xcolor}

\newmuskip\pFqmuskip

\newcommand*\pFq[6][8]{%
  \begingroup 
  \pFqmuskip=#1mu\relax
  \mathcode`=\string"8000
  \begingroup\lccode`\~=`\,
  \lowercase{\endgroup\let~}\pFqcomma
  F^{#2}_{#3}{\left(\genfrac..{0pt}{}{#4}{#5}\bigg|#6\right)}%
  \endgroup
}
\newcommand{\pFqcomma}{\mskip\pFqmuskip}

\newtheorem{theorem}{Theorem}
\newtheorem{lemma}[theorem]{Lemma}

\begin{document}

\title[Degenerate Binomial and degenerate Poisson Random Variables]{Degenerate Binomial and degenerate Poisson Random Variables}

\author{Dae San Kim}
\address{Department of Mathematics, Sogang University, Seoul 121-742, Republic of Korea}
\email{dskim@sogang.ac.kr}

\author{Taekyun  Kim}
\address{Department of Mathematics, Kwangwoon University, Seoul 139-701, Republic of Korea}
\email{tkkim@kw.ac.kr}


\maketitle

\begin{abstract}
The aim of this paper is to study the Poisson random variables in relation to the Lah-Bell polynomials and the degenerate binomial and degenerate Poisson random variables in connection with the degenerate Lah-Bell polynomials. Among other things, we show that the rising factorial moments of the degenerate Poisson random variable with parameter $\alpha$ are given by the degenerate Lah-Bell polynomials evaluated at $\alpha$. We also show that the probability-generating function of the degenerate Poisson random variable is equal to the generating function of the degenerate Lah-Bell polynomials. Also, we show similar results for the Poisson random variables. Here the $n$th Lah-Bell number counts the number of ways a set of $n$ elements can be partitioned into non-empty linearly ordered subsets, the Lah-Bell polynomials are natural extensions of the Lah-Bell numbers and the degenerate Lah-Bell polynomials are degenerate versions of the Lah-Bell polynomials.
\end{abstract}
 
\section{Introduction}

The aim of this paper is to study the Poisson random variables in relation to the Lah-Bell polynomials and the degenerate binomial and degenerate Poisson random variables in connection with the degenerate Lah-Bell polynomials. Here Lah-Bell polynomials $B_{n}^{L}(x)$ are natural extension of the Lah-Bell numbers $B_{n}^{L}$, which are defined as the number of ways a set of $n$ elements can be partitioned into non-empty linearly ordered subsets (see [3]). Thus we have $B_{n}^{L}=\sum_{k=0}^{n}L(n,k)$, where $L(n,k)$ counts the number of ways a set of $n$ elements can be partitioned into $k$ nonempty linearly ordered subsets, called the unsigned Lah numbers (see (2)). The motivation for our introduction of the degenerate binomial and degenerate Poisson random variables is as follows. Let us assume that the probability of success in a trial is $p$. Then we might wonder if the probability of success in the ninth trial is still $p$ after failing eight times in the trial experiment. Because there is a psychological burden for one to be successful. It seems plausible that the probability is less than $p$. This speculation motivated our study of the degenerate binomial and degenerate Poisson random variables. \par 
The outline of our main results is as follows. We derive the expectation and variance of the degenerate binomial and degenerate Poissson random variables. Then we introduce the degenerate Lah-Bell polynomials which are degenerate versions of the Lah-Bell polynomials. Then, among other things, we show that the rising factorial moments of the degenerate Poisson random variable with parameter $\alpha$ are given by the degenerate Lah-Bell polynomials evaluated at $\alpha$. We also show that the probability-generating function of the degenerate Poisson random variable is equal to the generating function of the degenerate Lah-Bell polynomials. In addition, we show that the rising factorial moments of the Poisson random variable with parameter $\alpha$ are given by the Lah-Bell polynomials evaluated at $\alpha$. Further, we show that the probability-generating function of the Poisson random variable is equal to the generating function of the Lah-Bell polynomials. \par
The novelty of this paper is that it reveals the connection between the rising factorial moments of the Poisson random variable and the Lah-Bell polynomials and that between the rising factorial moments of the degenerate Poisson random variable and the degenerate Lah-Bell polynomials. For the rest of this section, we recall the necessary facts that will be needed throughout this paper.

\vspace{0.1in}

For any $0 \ne \lambda\in\mathbb{R}$, the degenerate exponential functions are defined by 
\begin{equation}
e_{\lambda}^{x}(t)=\sum_{k=0}^{\infty}\frac{(x)_{k,\lambda}}{k!}t^{k},\quad(\mathrm{see}\ [4]), \label{1}
\end{equation}
where $(x)_{0,\lambda}=1,\ (x)_{n,\lambda}=x(x-\lambda)\cdots(x-(n-1)\lambda),\ (n\ge 1)$. 
Note that 
\begin{displaymath}
	\lim_{\lambda\rightarrow 0}e_{\lambda}^{x}(t)=e^{xt},\quad e_{\lambda}(t)=e_{\lambda}^{1}(t). 
\end{displaymath} \par
For $n,k\ge 0$, the unsigned Lah numbers are given by 
\begin{equation}
	L(n,k)=\binom{n-1}{k-1}\frac{n!}{k!},\quad(\mathrm{see}\ [2,3,7,10]). \label{2}
\end{equation}
In [3], the Lah-Bell polynomials are defined by 
\begin{equation}
e^{x\big(\frac{1}{1-t}-1\big)}=\sum_{n=0}^{\infty}B_{n}^{L}(x)\frac{t^{n}}{n!}. \label{3}
\end{equation}
For $x=1$, $B_{n}^{L}=B_{n}^{L}(1)$, $(n\ge 0)$, are called the Lah-Bell numbers. Here we recall from [3] that $B_{n}^{L}$ counts the number of ways a set of $n$ elements can be partitioned into non-empty linearly ordered subsets. 
From \eqref{3}, we note that 
\begin{equation}\label{4}
B_{n}^{L}(x)=\sum_{k=0}^{n}x^{k}L(n,k),\quad(n\ge 0),\quad(\mathrm{see}\ [3]). 
\end{equation}  \par 
A sample space is the set of all possible outcomes of an experiment and an event is any subset of the sample space. A random variable $X$ is a real valued function on a sample space. If $X$ takes any values in a countable set, then $X$ is called a discrete random variable. If $X$ takes any values in an interval on the real line, then $X$ is called a continuous random variable. \par 
For a discrete random variable $X$, the probability mass function $p(a)$ of $X$ is defined as 
\begin{equation}
p(a)=P\{X=a\},\quad(\mathrm{see}\ [8]). \label{5}
\end{equation}
Suppose that $n$ independent trials, each of which results in a ``success" with probability $p$ and in a ``failure" with probability $1-p$, are to be performed. If $X$ denotes the number of successes that occur in $n$ trials, then $X$ is called the binomial random variable with parameter $n,p$, which is denoted by $X\sim B(n,p)$. Let $X\sim B(n,p)$. Then the probability mass function of $X$ is given by 
\begin{equation}
	p(i)=\binom{n}{i}p^{i}(1-p)^{n-i},\quad i=0,1,2,\dots,n. \label{6}
\end{equation} \par
A Poisson random variable indicates how many events occured within a given period of time. A random variable $X$, taking on one of the values $0,1,2,\dots$, is said to be the Poisson random variable with parameter $\alpha(>0)$, if the probability mass function of $X$ is given by 
\begin{equation}
	p(i)=e^{-\alpha}\frac{\alpha^{i}}{i!},\quad(\mathrm{see}\ [8]). \label{7}
\end{equation} \par
Let $f(x)$ be a real valued function and let $X$ be a random variable. Then we define 
\begin{equation}
	E[f(X)]=\sum_{i=0}^{\infty}f(i)p(i),\quad(\mathrm{see}\ [8]), \label{8}
\end{equation}
where $p(x)$ is the probability mass function of $X$. \par 
It is well known that the Bell polynomials are defined by 
\begin{equation}
	e^{x(e^{t}-1)}=\sum_{n=0}^{\infty}B_{n}(x)\frac{t^{n}}{n!},\quad(\mathrm{see}\ [1,4,5]). \label{9}
\end{equation} 
Let us take $f(x)=x^{n},\ (n\ge 0)$. Then we have the moments of the Poisson random variable $X$ with parameter $\alpha(>0)$ as follows : 
\begin{equation}
	E[X^{n}]=B_{n}(\alpha),\quad(n\ge 0),\quad(\mathrm{see}\ [5]). \label{10}
\end{equation}

\section{Poisson random variables}

The falling factorial sequence is given by 
\begin{displaymath}
(x)_{0}=1,\quad (x)_{n}=x(x-1)\cdots(x-n+1),\quad(n\ge 1), 
\end{displaymath}
while the rising factorial sequence is given by 
\begin{displaymath}
\langle x\rangle_{0}=1,\quad \langle x\rangle_{n}=x(x+1)\cdots(x+n-1),\quad(n\ge 1),\quad(\mathrm{see}\ [1-10]).
\end{displaymath}
Replacing $t$ by $\log(1+t)$ in \eqref{9}, we get 
\begin{align}
e^{xt}\ &=\ \sum_{k=0}^{\infty}B_{k}(x)\frac{1}{k!}\big(\log(1+t)\big)^{k} \label{11} \\
&=\ \sum_{k=0}^{\infty}B_{k}(x)\sum_{n=k}^{\infty}S_{1}(n,k)\frac{t^{n}}{n!} \nonumber \\
&=\ \sum_{n=0}^{\infty}\bigg(\sum_{k=0}^{n}B_{k}(x)S_{1}(n,k)\bigg)\frac{t^{n}}{n!}, \nonumber
\end{align}
where $S_{1}(n,k)$ are the Stirling numbers of the first kind defined by 
\begin{equation}
	(x)_{n}=\sum_{k=0}^{n}S_{1}(n,k)x^{k},\quad(n\ge 0). \label{12}
\end{equation}
Therefore, by \eqref{11}, we obtain the following lemma.
\begin{lemma}
For $n\ge 0$, we have 
\begin{displaymath}
x^{n}=\sum_{k=0}^{n}S_{1}(n,k)B_{k}(x), 
\end{displaymath}
and 
\begin{displaymath}
B_{n}(x)=\sum_{k=0}^{n}S_{2}(n,k)x^{k}. 
\end{displaymath}
\end{lemma}
Let $X$ be the Poisson random variable with parameter $\alpha(>0)$. Then we have 
\begin{align}
E[(X)_{n}]\ &=\ \sum_{k=0}^{n}S_{1}(n,k)E[X^{k}] \nonumber \\
&=\ \sum_{k=0}^{n}S_{1}(n,k)B_{k}(\alpha). \label{13}
\end{align}
From Lemma 1 and \eqref{13}, we note the well-known fact about the falling factorial moments of the random variable $X$, namely the expectation of the falling factorial of the random variable $X$:
\begin{equation}
E[(X)_{n}]=\sum_{k=0}^{n}S_{1}(n,k)B_{k}(\alpha)=\alpha^{n},\quad(n\ge 0). \label{14}
\end{equation} \par
On the other hand, the rising factorial moment of $X$, namely the expectation of the rising factorial of $X$, is given by
\begin{align}
E[\langle X\rangle_{n}]\ &=\ \sum_{k=0}^{\infty}\langle k\rangle_{n}p(k) \label{15} \\
&=\ e^{-\alpha}\sum_{k=0}^{\infty}\frac{\langle k\rangle_{n}}{k!}\alpha^{k}. \nonumber
\end{align}
From \eqref{3}, we can derive the following equation. 
\begin{align}
\sum_{n=0}^{\infty}B_{n}^{L}(\alpha)\frac{t^{n}}{n!}\ &=\ e^{-\alpha}e^{\alpha\big(\frac{1}{1-t}\big)} \label{16} \\
&=e^{-\alpha}\sum_{k=0}^{\infty}\alpha^{k}\frac{1}{k!}\bigg(\frac{1}{1-t}\bigg)^{k} \nonumber \\
&=\ e^{-\alpha}\sum_{k=0}^{\infty}\frac{\alpha^{k}}{k!}\sum_{n=0}^{\infty}\langle k\rangle_{n}\frac{t^{n}}{n!} \nonumber \\
&=\ \sum_{n=0}^{\infty}\bigg(e^{-\alpha}\sum_{k=0}^{\infty}\frac{\langle k\rangle_{n}}{k!}\alpha^{k}\bigg)\frac{t^{n}}{n!}. \nonumber
\end{align}
Comparing the coefficients on both sides of \eqref{16}, we have the following identity:
\begin{equation}
B_{n}^{L}(\alpha)=e^{-\alpha}\sum_{k=0}^{\infty}\frac{\langle k\rangle_{n}}{k!}\alpha^{k},\label{17}
\end{equation}
where $n$ is a nonnegative integer. \par 
Therefore, by \eqref{14}, \eqref{15} and \eqref{17}, we obtain the following theorem. In particular, it shows that the rising factorial moments of the Poisson random variable with parameter $\alpha$ are given by the Lah-Bell polynomials evaluated at $\alpha$. This fact seems to be new.
\begin{theorem}
Let $X$ be the Poisson random variable with parameter $\alpha(>0)$. Then we have
\begin{displaymath}
E[(X)_{n}]=\alpha^{n}, 
\end{displaymath}
and 
\begin{displaymath}
E[\langle X\rangle_{n}]=B_{n}^{L}(\alpha), \quad(n\ge 0). 
\end{displaymath}
\end{theorem}
Let $X$ be the Poission random variable with parameter $\alpha(>0)$. From \eqref{7} and \eqref{8}, we have 
\begin{align}
E\bigg[\bigg(\frac{1}{1-t}\bigg)^{X}\bigg]\ &=\ \sum_{k=0}^{\infty}\bigg(\frac{1}{1-t}\bigg)^{k}p(k) \nonumber \\
&=\ \sum_{k=0}^{\infty}\bigg(\frac{1}{1-t}\bigg)^{k}e^{-\alpha}\frac{\alpha^{k}}{k!} \label{18} \\
&=\ e^{-\alpha}e^{\frac{\alpha}{1-t}}\ =\ e^{\alpha\big(\frac{1}{1-t}-1\big)}.\nonumber
\end{align}
Now, by \eqref{3} and \eqref{18}, we obtain the following theorem which says that the probability-generating function of $X$ is equal to the generationg function of the Lah-Bell polynomials.
\begin{theorem}
Let $X$ be the Poission random variable with parameter $\alpha(>0)$. Then we have 
\begin{displaymath}
E\bigg[\bigg(\frac{1}{1-t}\bigg)^{X}\bigg]\ =\ e^{\alpha\big(\frac{1}{1-t}-1\big)}\ =\ \sum_{n=0}^{\infty}B_{n}^{L}(\alpha)\frac{t^{n}}{n!}. 
\end{displaymath}
\end{theorem}
From Theorem 3 and \eqref{10}, we note that 
\begin{align}
\sum_{n=0}^{\infty}B_{n}^{L}(\alpha)\frac{t^{n}}{n!}\ & =\ \sum_{k=0}^{\infty}E[X^{k}]\frac{(-\log(1-t))^k}{k!} \label{19}\\
&=\ \sum_{k=0}^{\infty}(-1)^{k}B_{k}(\alpha)\sum_{n=k}^{\infty}(-1)^{n}S_{1}(n,k)\frac{t^n}{n!}\nonumber \\
&=\ \sum_{n=0}^{\infty}\bigg(\sum_{k=0}^{n}(-1)^{n-k}S_{1}(n,k)B_{k}(\alpha)\bigg)\frac{t^n}{n!}.\nonumber
\end{align}
Therefore, by Theorem 2, \eqref{4} and \eqref{19}, we obtain the following theorem. 
\begin{theorem}
Let $X$ be the Poission random vairable with parameter $\alpha(>0)$. Then we have 
\begin{displaymath}
E[\langle X\rangle_{n}] \ =\ B_{n}^{L}(\alpha)\ =\ \sum_{k=0}^{n}L(n,k)\alpha^{k}\ =\ \sum_{k=0}^{n}(-1)^{n-k}S_{1}(n,k)B_{k}(\alpha). 
\end{displaymath}
\end{theorem}

\section{Degenerate Binomial and degenerate Poisson Random Variables} 
In this section, we assume that $\lambda\in (0,1)$, and $p$ is the probability of success in an experiment. 
For $\lambda\in (0,1)$, $X_{\lambda}$ is the {\it{degenerate binomial random variable with parameter $n,p$}}, denoted by $X_{\lambda}\sim B_{\lambda}(n,p)$, if the probability mass function of $X_{\lambda}$ is given by 
\begin{equation}
p_{\lambda}(i)=P\{X_{\lambda}=i\}=\binom{n}{i}(p)_{i,\lambda}(1-p)_{n-i,\lambda}\frac{1}{(1)_{n,\lambda}},\label{20}
\end{equation}
where $i=0,1,2,\dots,n$. \par 
From \eqref{20}, we note that 
\begin{displaymath}
	\sum_{i=0}^{\infty}p_{\lambda}(i)=\frac{1}{(1)_{n,\lambda}}\sum_{i=0}^{n}\binom{n}{i}(p)_{i,\lambda}(1-p)_{n-i,\lambda}=1. 
\end{displaymath}
For $X_{\lambda}\sim B_{\lambda}(n,k)$, we have 
\begin{align}
E[X_{\lambda}]\ &=\ \sum_{i=0}^{\infty}ip_{\lambda}(i) \label{21}\\ 
&=\ \frac{1}{(1)_{n,\lambda}}\sum_{i=0}^{\infty}i\binom{n}{i}(p)_{i,\lambda}(1-p)_{n-i,\lambda}\nonumber \\
&=\ \frac{n}{(1)_{n,\lambda}}\sum_{i=1}^{\infty}\binom{n-1}{i-1}(p)_{i,\lambda}(1-p)_{n-i,\lambda}\nonumber\\
&=\ \frac{n}{(1)_{n,\lambda}}\sum_{i=0}^{\infty}\binom{n-1}{i}(p)_{i+1,\lambda}(1-p)_{n-1-i,\lambda} \nonumber\\
&=\ \frac{np}{(1)_{n,\lambda}}\sum_{i=0}^{n-1}\binom{n-1}{i}(p-\lambda)_{i,\lambda}(1-p)_{n-1-i,\lambda} \nonumber \\
&=\ \frac{np}{(1)_{n,\lambda}}(p-\lambda+1-p)_{n-1,\lambda}\nonumber \\
&= \frac{np}{(1)_{n,\lambda}}(1-\lambda)_{n-1,\lambda}. \nonumber
\end{align}
Therefore, we obtain the following theorem. 
\begin{theorem}
For $X_{\lambda}\sim B_{\lambda}(n,p),\ (n\ge 0)$, we have 
\begin{displaymath}
E[X_{\lambda}]=\frac{np}{(1)_{n,\lambda}}(1-\lambda)_{n-1,\lambda}. 
\end{displaymath}
\end{theorem}
Note that 
\begin{displaymath}
\lim_{\lambda\rightarrow 0}E[X_{\lambda}]=np=E[X], 
\end{displaymath}
where $X$ is the binomial random variable with parameter $n,p$. \par 
For $X_{\lambda}\sim B_{\lambda}(n,p)$, we observe that 
\begin{align}
E[X_{\lambda}^{2}]\ &=\ \sum_{i=0}^{\infty}i^{2}p_{\lambda}(i)\ =\ \frac{1}{(1)_{n,\lambda}}\sum_{i=0}^{\infty}i^{2}\binom{n}{i}(p)_{i,\lambda}(1-p)_{n-i,\lambda} \label{22} \\
&=\ \frac{1}{(1)_{n,\lambda}}\sum_{i=0}^{\infty}i(i-1+1)\binom{n}{i}(p)_{i,\lambda}(1-p)_{n-i,\lambda}\nonumber \\
&=\ \frac{1}{(1)_{n,\lambda}}\sum_{i=0}^{\infty}i(i-1)\binom{n}{i}(p)_{i,\lambda}(1-p)_{n-i,\lambda}+E[X_{\lambda}]\nonumber  \\
&=\ \frac{n(n-1)}{(1)_{n,\lambda}}\sum_{i=2}^{\infty}\binom{n-2}{i-2}(p)_{i,\lambda}(1-p)_{n-i,\lambda}+E[X_{\lambda}] \nonumber \\
&=\ \frac{n(n-1)}{(1)_{n,\lambda}}\sum_{i=0}^{n-2}\binom{n-2}{i}(p)_{i+2,\lambda}(1-p)_{n-2-i,\lambda}+E[X_{\lambda}] \nonumber \\
&=\ \frac{n(n-1)p(p-\lambda)}{(1)_{n,\lambda}}\sum_{i=0}^{n-2}\binom{n-2}{i}(p-2\lambda)_{i,\lambda}(1-p)_{n-2-i,\lambda}+E[X_{\lambda}] \nonumber \\
&=\ \frac{n(n-1)p(p-\lambda)}{(1)_{n,\lambda}}\big(p-2\lambda+1-p\big)_{n-2,\lambda}+E[X_{\lambda}] \nonumber \\
&=\ \frac{n(n-1)p(p-\lambda)}{(1)_{n,\lambda}}(1-2\lambda)_{n-2,\lambda}+\frac{np}{(1)_{n,\lambda}}(1-\lambda)_{n-1,\lambda} \nonumber \\
&=\ \frac{np(1-2\lambda)_{n-2,\lambda}}{(1)_{n,\lambda}}\big\{(p-\lambda)(n-1)+(1-\lambda)\big\}\nonumber \\
&=\ \frac{np}{(1)_{n,\lambda}}(1-2\lambda)_{n-2,,\lambda}\big(p(n-1)+1-n\lambda\big). \nonumber 
\end{align}

By using Theorem 5 and \eqref{22}, the variance $\mathrm{Var}(X_{\lambda})$ of the random variable $X_{\lambda}$ is given by
\begin{align}
\mathrm{Var}(X_{\lambda})&=E[X_{\lambda}^{2}]-\big(E[X_{\lambda}]\big)^{2} \label{23}\\
&=\ \frac{np}{(1)_{n,\lambda}}(1-2\lambda)_{n-2,\lambda}\big(p(n-1)+1-n\lambda)-\bigg(\frac{np}{(1)_{n,\lambda}}(1-\lambda)_{n-1,\lambda}\bigg)^{2}\nonumber\\
&=\ \frac{np}{(1)_{n,\lambda}}(1-2\lambda)_{n-2,\lambda}\big(p(n-1)+1-n\lambda\big)-\bigg(\frac{np}{(1)_{n,\lambda}}(1-\lambda)(1-2\lambda)_{n-2,\lambda}\bigg)^{2}\nonumber\\
&=\ \frac{np(1-2\lambda)_{n-2,\lambda}}{(1)_{n,\lambda}}\bigg(p(n-1)+1-n\lambda-\frac{np}{(1)_{n,\lambda}}(1-\lambda)^{2}(1-2\lambda)_{n-2,\lambda}\bigg)\nonumber\\
&=\ \frac{np}{(1)_{n,\lambda}}(1-2\lambda)_{n-2,\lambda}\big((n-1)p+1-n\lambda-E[X_{\lambda}](1-\lambda)\big).\nonumber
\end{align}
Therefore, we obtain the following theorem. 
\begin{theorem}
For $X_{\lambda}\sim B_{\lambda}(n,p)$, we have 
\begin{displaymath}
\mathrm{Var}(X_{\lambda})=\frac{np}{(1)_{n,\lambda}}(1-2\lambda)_{n-2,\lambda}\big((n-1)p+1-n\lambda-E[X_{\lambda}](1-\lambda)\big). 
\end{displaymath}
\end{theorem}
Note that 
\begin{displaymath}
	\lim_{\lambda\rightarrow 0}\mathrm{Var}(X_{\lambda})=np(1-p)=\mathrm{Var}(X), 
\end{displaymath}
where $X$ is the binomial random variable with parameters $n,p$. \par 
The generating function of the moments of $X_{\lambda} \sim B_{\lambda}(n,p)$ is given by 
\begin{align*}
\sum_{n=0}^{\infty}E[X_{\lambda}^{n}]\frac{t^{n}}{n!}\ &=\ E[e^{X_{\lambda} t}] \\
	&=\ \frac{1}{(1)_{n,\lambda}}\sum_{i=0}^{n}e^{it}\binom{n}{i}(p)_{i,\lambda}(1-p)_{n-i,\lambda}. 
\end{align*}
Thus, we have 
\begin{align*}
\ E[X_{\lambda}^{n}]&= \frac{d^{n}}{dt^{n}}E\big[e^{X_{\lambda}t}\big]\bigg|_{t=0} \\
&= \frac{1}{(1)_{n,\lambda}}\sum_{i=0}^{n}\binom{n}{i}i^{n}(p)_{i,\lambda}(1-p)_{n-i,\lambda}. 
\end{align*}
For $\lambda\in (0,1)$, $X_{\lambda}$ is the {\it{ degenerate Poisson random variable with parameter $\alpha(>0)$}}, if the probability mass function of $X_{\lambda}$ is given by 
\begin{equation}
p_{\lambda}(i)=P\{X_{\lambda}=i\}=e_{\lambda}^{-1}(\alpha)\frac{\alpha^{i}}{i!}(1)_{i,\lambda}, \label{24}
\end{equation}
where $i=0,1,2,\dots$. \par 
By \eqref{24}, we get 
\begin{displaymath}
\sum_{i=0}^{\infty}p_{\lambda}(i)=e_{\lambda}^{-1}(\alpha)\sum_{i=0}^{\infty}\frac{(1)_{i,\lambda}}{i!}\alpha^{i}=e_{\lambda}^{-1}(\alpha)e_{\lambda}(\alpha)=1. 
\end{displaymath}
It is easy to show that 
\begin{displaymath}
E[X_{\lambda}]=\frac{\alpha}{1+\alpha\lambda},
\end{displaymath}
and 
\begin{displaymath}
E[X_{\lambda}^{2}]=\frac{\alpha+\alpha^{2}}{(1+\alpha\lambda)^{2}}. 
\end{displaymath}
Thus, we have 
\begin{displaymath}
\mathrm{Var}(X_{\lambda})=E[X_{\lambda}^{2}]-\big(E[X_{\lambda}]\big)^{2}=\frac{\alpha}{(1+\alpha\lambda)^{2}}. 
\end{displaymath} \par
Let $X_{\lambda}$ be the degenerate Poisson random variable with parameter $\alpha(>0)$. 
Then we have
\begin{align}
E[\langle X_{\lambda}\rangle_{n}]\ &=\ \sum_{i=0}^{\infty}\langle i\rangle _{n}p_{\lambda}(i) \label{25} \\
&=\ \sum_{i=0}^{\infty}\langle i\rangle_{n}e_{\lambda}^{-1}(\alpha)\frac{(1)_{n,\lambda}}{i!}\alpha^{i}\nonumber \\
&= e^{-1}_{\lambda}(\alpha)\sum_{i=0}^{\infty}(1)_{i,\lambda}\frac{\langle i\rangle_{n}}{i!}\alpha^{i}. \nonumber 
\end{align} \par
In view of \eqref{3}, we may consider the {\it{degenerate Lah-Bell polynomials}} given by
\begin{equation}
e^{-1}_{\lambda}(x)e_{\lambda}\bigg(x\bigg(\frac{1}{1-t}\bigg)\bigg)=\sum_{n=0}^{\infty}B_{n,\lambda}^{L}(x)\frac{t^{n}}{n!}. \label{26}
\end{equation}
Note that 
\begin{displaymath}
\sum_{n=0}^{\infty}\lim_{\lambda\rightarrow 0}B_{n,\lambda}^{L}(x)\frac{t^{n}}{n!}=e^{x\big(\frac{1}{1-t}-1\big)}=\sum_{n=0}^{\infty}B_{n}^{L}(x)\frac{t^{n}}{n!}. 
\end{displaymath}
Thus we have
\begin{displaymath}
\lim_{\lambda\rightarrow 0}B_{n,\lambda}^{L}(x)=B_{n}^{L}(x),\quad(n\ge 0). 
\end{displaymath}
Now, we observe that 
\begin{align}
e_{\lambda}^{-1}(x)e_{\lambda}\bigg(x\bigg(\frac{1}{1-t}\bigg)\bigg)\ &=\ e_{\lambda}^{-1}(x)\sum_{k=0}^{\infty}(1)_{k,\lambda}x^{k}\frac{1}{k!}\bigg(\frac{1}{1-t}\bigg)^{k} \label{27} \\
&=\ e_{\lambda}^{-1}(x)\sum_{k=0}^{\infty}(1)_{k,\lambda}x^{k}\frac{1}{k!}\sum_{n=0}^{\infty}\frac{\langle k\rangle_{n}}{n!}t^{n} \nonumber \\
&=\ \sum_{n=0}^{\infty}\bigg\{e_{\lambda}^{-1}(x)\sum_{k=0}^{\infty}(1)_{k,\lambda}\frac{\langle k\rangle_{n}}{k!}x^{k}\bigg\}\frac{t^{n}}{n!}. \nonumber	
\end{align}
From \eqref{25}, \eqref{26} and \eqref{27}, we obtain the next result. In particular, it says that the rising factorial moments of the discrete Poisson random variable with parameter $\alpha$ are given by 
the degenerate Lah-Bell polynomials evaluated at $\alpha$. 
\begin{theorem}
For $n\ge 0$, we have 
\begin{displaymath}
B_{n,\lambda}^{L}(x)=e_{\lambda}^{-1}(x)\sum_{k=0}^{\infty}(1)_{k,\lambda}\frac{\langle k\rangle_{n}}{k!}x^{k}.
\end{displaymath}
In particular, we have
\begin{displaymath}
E[\langle X_{\lambda}\rangle_{n}]=B_{n,\lambda}^{L}(\alpha),\quad(n\ge 0). 
\end{displaymath}
\end{theorem}
The degenerate Bell polynomials are defined in [4] as 
\begin{equation}
e_{\lambda}^{-1}(x)e_{\lambda}(xe^{t})=\sum_{n=0}^{\infty}B_{n,\lambda}(x)\frac{t^{n}}{n!}. \label{28}
\end{equation}
Note that 
\begin{displaymath}
\sum_{n=0}^{\infty}\lim_{\lambda\rightarrow 0}B_{n,\lambda}(x)\frac{t^{n}}{n!}=e^{x(e^{t}-1)}=\sum_{n=0}^{\infty}B_{n}(x)\frac{t^{n}}{n!}, 
\end{displaymath}
where $B_{n}(x)$ are the ordinary Bell polynomials.  
Thus, we have 
\begin{displaymath}
\lim_{\lambda\rightarrow 0}B_{n,\lambda}(x)=B_{n}(x),\quad(n\ge 0). 
\end{displaymath} \par
Replacing $t$ by $-\log(1-t)$ in \eqref{28}, we get 
\begin{align}
e^{-1}_{\lambda}(x)e_{\lambda}\bigg(x\bigg(\frac{1}{1-t}\bigg)\bigg)\ &=\ \sum_{k=0}^{\infty}B_{k,\lambda}(x)(-1)^{k}\frac{1}{k!}\big(\log(1-t)\big)^{k} \nonumber\\
&=\ \sum_{k=0}^{\infty}B_{k,\lambda}(x)(-1)^{k}\sum_{n=k}^{\infty}(-1)^{n}S_{1}(n,k)\frac{t^n}{n!}\label{29}\\
&=\ \sum_{n=0}^{\infty}\bigg(\sum_{k=0}^{n}(-1)^{n-k}S_{1}(n,k)B_{k,\lambda}(x)\bigg)\frac{t^{n}}{n!}. \nonumber
\end{align}
Therefore, by \eqref{26}, \eqref{28} and \eqref{29}, we obtain the following theorem. 
\begin{theorem}
For $n\ge 0$, we have 
\begin{equation}\label{30}
B_{n,\lambda}^{L}(x)=\sum_{k=0}^{n}(-1)^{n-k}S_{1}(n,k)B_{k,\lambda}(x),
\end{equation}
and 
\begin{displaymath}
B_{n,\lambda}(x)=\sum_{k=0}^{n}(-1)^{n-k}S_{2}(n,k)B_{k,\lambda}^{L}(x), 
\end{displaymath}
where $S_{2}(n,k)$, $(n,k\ge 0)$, are the Stirling numbers of the second kind defined by 
\begin{displaymath}
x^{n}=\sum_{k=0}^{n}S_{2}(n,k)(x)_{k}.
\end{displaymath}
\end{theorem}
From Theorem 11 of [4], we recall that 
\begin{equation}\label{31}
B_{n,\lambda}(x)=\sum_{k=0}^{n}(1)_{k,\lambda}\bigg(\frac{x}{1+\lambda x}\bigg)^{k}S_{2}(n,k).
\end{equation}
Combining \eqref{30} and \eqref{31}, we have another expression for $B_{n,\lambda}(x)$ as follows.
\begin{equation}\label{32}
B_{n,\lambda}^{L}(x)=\sum_{l=0}^{n}\bigg(\sum_{k=l}^{n}(-1)^{n-k}S_1(n,k)S_2(k,l)\bigg)(1)_{l,\lambda}\bigg(\frac{x}{1+\lambda x}\bigg)^l.
\end{equation}
Let $X_{\lambda}$ be the degenerate Poisson random variable with parameter $\alpha(>0)$. Then we have
\begin{align}
\ E\bigg[\bigg(\frac{1}{1-t}\bigg)^{X_{\lambda}}\bigg]\ & =\ \sum_{i=0}^{\infty}\bigg(\frac{1}{1-t}\bigg)^{i}p_{\lambda}(i)\label{33}\\
&=\ e_{\lambda}^{-1}(\alpha)\sum_{i=0}^{\infty}(1)_{i,\lambda}\frac{\alpha^{i}}{i!}\bigg(\frac{1}{1-t}\bigg)^{i}\nonumber \\
&=\ e^{-1}_{\lambda}(\alpha)e_{\lambda}\bigg(\alpha\bigg(\frac{1}{1-t}\bigg)\bigg). \nonumber
\end{align}
Therefore, we obtain the following theorem from Theorem 7, \eqref{32} and \eqref{33}. In particular, it states that the probability-generating function of $X_{\lambda}$ is equal to the generating function of the degenerate Lah-Bell polynomials.
\begin{theorem}
Let $X_{\lambda}$ be the Poisson random variable with parameter $\alpha>0$. Then we have 
\begin{displaymath}
E\bigg[\bigg(\frac{1}{1-t}\bigg)^{X_{\lambda}}\bigg]=e_{\lambda}^{-1}(\alpha)\cdot e_{\lambda}\bigg(\alpha\bigg(\frac{1}{1-t}\bigg)\bigg),
\end{displaymath}
and
\begin{displaymath}
E[\langle X_{\lambda}\rangle_{n}]=\sum_{l=0}^{n}\bigg(\sum_{k=l}^{n}(-1)^{n-k}S_1(n,k)S_2(k,l)\bigg)(1)_{l,\lambda}\bigg(\frac{\alpha}{1+\lambda \alpha}\bigg)^l,\quad(n\ge 0). 
\end{displaymath}
\end{theorem}

\section{Conclusion}

In this paper, we derived the expectation and variance of the degenerate binomial and degenerate Poissson random variables. Then we introduced the degenerate Lah-Bell polynomials which are degenerate versions of the recently introduced Lah-Bell polynomials (see [3]). Then we showed that the rising factorial moments of the degenerate Poisson random variable with parameter $\alpha$ are given by the degenerate Lah-Bell polynomials evaluated at $\alpha$. We also showed that the probability-generating function of the degenerate Poisson random variable is equal to the generating function of the degenerate Lah-Bell polynomials. We also derived similar results for the Poisson random variable. \par
Here we would like to mention that studying various degenerate versions of some special numbers of polynomials, which was initiated by Carlitz when he investigated the degenerate Bernoulli and Euler polynomials and numbers, regained interests of some mathematicians in recent years. They have been studied by using several different tools like generating functions, combinatorial methods, $p$-adic analysis, umbral calculus, special functions, differential equations and probability theory as we did in the present paper.
It is one of our future projects to continue to study various degenerate versions of some special polynomials and numbers and to find their applications in physics, science and engineering as well as in mathematics.

\end{document}